\documentclass{amsart}

\usepackage{amssymb,amsfonts,amsxtra,amsthm}
\usepackage{dsfont,mathrsfs}
\renewcommand{\mathcal}{\mathscr}
\usepackage[all]{xypic}
\xyoption{dvips}
\usepackage{epic,eepic}
\usepackage{paralist,longtable}

\newcommand{\CC}{\mathds{C}}
\newcommand{\Cs}{\CC\{\!\{\del_t^{-1}\}\!\}}
\newcommand{\calD}{\mathcal{D}}
\newcommand{\calG}{\mathcal{G}}
\newcommand{\calH}{\mathcal{H}}
\newcommand{\calO}{\mathcal{O}}
\newcommand{\gl}{\mathfrak{gl}}
\newcommand{\del}{\partial}
\newcommand{\HH}{\mathbf{H}}
\newcommand{\ideal}[1]{{\langle#1\rangle}}
\newcommand{\into}{\hookrightarrow}
\newcommand{\Lie}{\mathrm{L}}

\newcommand{\mm}{\mathfrak{m}}
\newcommand{\NN}{\mathds{N}}
\newcommand{\QQ}{\mathds{Q}}
\newcommand{\RR}{\mathbf{R}}
\newcommand{\xymat}{\SelectTips{cm}{}\xymatrix}

\DeclareMathOperator{\Der}{Der}
\DeclareMathOperator{\gr}{gr}

\DeclareMathOperator{\tr}{tr}

\theoremstyle{definition}
\newtheorem{ntn}{Notation}

\theoremstyle{plain}

\newtheorem{prp}[ntn]{Proposition}
\newtheorem{thm}[ntn]{Theorem}
\newtheorem{cor}[ntn]{Corollary}

\theoremstyle{remark}

\newtheorem{exa}[ntn]{Example}

\begin{document}

\title[Logarithmic comparison for isolated singularities]{Logarithmic comparison theorem versus Gauss--Manin system for isolated singularities}

\author{Mathias Schulze}

\address{M.~Schulze\\
Oklahoma State University\\
Department of Mathematics\\
Stillwater, OK 74078\\
United States}
\email{mschulze@math.okstate.edu}

\thanks{The author is grateful to Claus Hertling and to the anonymous referee for careful reading and useful suggestions, and to Michel Granger for helpful discussions.}

\date{\today}

\begin{abstract}
For quasihomogeneous isolated hypersurface singularities, the logarithmic comparison theorem has been characterized explicitly by Holland and Mond.
In the nonquasihomogeneous case, we give a necessary condition for the logarithmic comparison theorem in terms of the Gauss--Manin system of the singularity.
It shows in particular that the logarithmic comparison theorem can hold for a nonquasihomogeneous singularity only if $1$ is an eigenvalue of the monodromy.
\end{abstract}


\maketitle

\section{Introduction}

Let $D$ be a hypersurface in a complex manifold $X$ with complement $j:U=X\setminus D\into X$.
Then Grothendieck's comparison theorem \cite{Gro66} states that the De~Rham morphism
\[
\Omega^\bullet_X(*D)\to\RR j_*\CC_U
\]
is a quasi--isomorphism.
In particular, for Stein $X$, each cohomology class $c\in H^k(U;\CC)$ is represented as $c(\sigma)=\int_\sigma\omega$ by a differential $k$-form $\omega$ with finite pole order along $D$.
The natural question of limiting this pole order dates back to Griffiths \cite{Gri69} and has been studied later by Deligne and Dimca \cite{DD90,Dim91}, Karpishpan \cite{Kar91}, and others.

For a normal crossing divisor $D=\{x_1\cdots x_k=0\}\subseteq\CC^n=X$, the poles can be restricted to simple poles. 
More precisely, the inclusion 
\begin{equation}\label{7}
\Omega^\bullet(\log D)\into\Omega^\bullet_X(*D)
\end{equation}
of the complex of logarithmic differential forms, generated by $\frac{dx_1}{x_1},\dots,\frac{dx_k}{x_k},dx_{k+1},\dots,dx_n$, is a quasi--isomorphism.
This fact plays a crucial role in Deligne's mixed Hodge theory \cite[\S 3]{Del71}.

Saito \cite{Sai80} extended the definition of the complex $\Omega^\bullet(\log D)$ of logarithmic differential forms to general hypersurfaces $D$.
In analogy with Grothendieck's theorem, one says that the logarithmic comparison theorem holds for $D$ if \eqref{7} is a quasi--isomorphism.
The problem of characterizing such $D$ has been studied essentially in the extremal cases of isolated singularities and of free divisors.
In both cases the complete solution is still missing.

The overlap of the two cases, the plane curve case, is completely understood: The logarithmic comparison theorem is equivalent to quasihomogeneity of the singularities \cite{CMNC02}.
For free divisors, the normal crossing case has been extended to the class of (weakly) locally quasihomogeneous free divisors, for which the logarithmic comparison theorem holds \cite{CNM96} (\cite[Rem.~1.7.4]{Nar08}).
For general free divisors, there is a D-module theoretic reformulation of the logarithmic comparison theorem based on a $\calD_X(-\log D)$-duality analogous to the ordinary $\calD_X$-duality \cite{CN05}.

In the present note we are concerned with the case of isolated singularities.
By the local nature of the problem, we can reduce to germs of spaces an maps:
\begin{gather}\label{22}
\xymat{X:=(\CC^{n+1},0)\ar[r]^-f&(\CC,0)=:T},\\
\nonumber f\in\mm^2\subseteq\mm:=\mm_X=\ideal{x}\subseteq\calO:=\calO_X=\CC\{x\},\quad\calO_T=\CC\{t\},\\
\nonumber D=\{x\mid f(x)=0\}\subseteq X,\quad\xymat{U:=X\setminus D\ar[r]^-j&X},
\end{gather}
where $f$ is a reduced equation of the isolated hypersurface singularity $D$, $x=x_0,\dots,x_n$ and $t$ are coordinates on $X$ and $T$.
We shall tacitly identify $X$ with a Milnor representative \cite{Mil68}.
Note that the latter form a basis of Stein neighborhoods of $0\in X$ and it suffices to check the logarithmic comparison theorem on global sections over such neighborhoods by \cite[Lem.~2.5]{CNM96}.

The main result for isolated singularities due to Holland and Mond \cite{HM98} covers the case of quasihomogeneous singularities.

\begin{thm}[Holland, Mond]\label{19}
Let $D=\{x\mid f(x)=0\}\subseteq X$ be a quasihomogeneous isolated hypersurface singularity of degree $r$ with respect to positive integer weights $w_1,\dots,w_n$.
Denote by $J_f$ the gradient ideal of $f$.
Then the following conditions are equivalent:
\begin{enumerate}[(a)]
\item the logarithmic comparison theorem holds at $0$;
\item $(\RR^ij_*\CC_U)_0=0$ for $i\ge2$;
\item $(\calO_{X,0}/J_f)_{ir-\sum_jw_j}=0$ for $1\le i\le n-1$;
\item the link of $0$ in $D$ is a $\QQ$-homology sphere.
\end{enumerate}
Furthermore, each of these statements implies that $\calH^i(\Omega^\bullet(\log D))=0$ for $i\ge2$; for $n=2$ the reverse implication also holds.
\end{thm}

For free divisors, it is conjectured, and proved for $n\le2$, that the logarithmic comparison theorem requires strong Euler homogeneity \cite{GS06}.
For isolated singularities, the latter property reduces to quasihomogeneity by \cite{Sai71} and one could expect that the logarithmic comparison theorem requires quasihomogeneity.
Our main result confirms this expectation for a large class of isolated singularities defined by properties of the Gauss--Manin system $\calG:=\int_f^0\calO_X$.
This is the direct image of the $\calD_X$-module $\calO_X$ along $f\colon X\to T$ and as such a $\calD_T$-module.

\begin{thm}\label{14}
Let $D=\{x\mid f(x)=0\}\subseteq X$ be an isolated hypersurface singularity with Gauss--Manin system $\calG=\int_f^0\calO_X$, Brieskorn lattice $\calH''\subseteq\calG$, monodromy $M$, and spectrum $\alpha_1<\alpha_2\le\cdots\le\alpha_n$.
Denote by $V$ the Kashiwara--Malgrange filtration on $\calG$, by $C^\alpha\cong\gr_V\calG$ the generalized $\alpha$-eigenspace of $t\del_t$ in $\calG$, by $H_\alpha$ the image of $\calH''\cap V^\alpha$ in $C^\alpha$, and $N=\log(M_u)\colon C^\alpha\to C^\alpha$ where $M_u$ is the unipotent part of $M$.

Under each of the following conditions the logarithmic comparison theorem can hold for $D$ only if $D$ is quasihomogeneous. 
\begin{enumerate}[(a)]
\item\label{14a} $1$ is not an eigenvalue of $M$;
\item\label{14b} $\alpha_1>0$;
\item\label{14c} $\alpha_1<0$ and $[udx]_0\in H_0\oplus N(C^0)$ for some $u\in\calO_X^*$;
\item\label{14d} $\alpha_1<0=\alpha_2$ and $[dx]_0\in H_0\oplus N(C^0)$.
\end{enumerate}
\end{thm}

In the case $\alpha_1=0$, our approach does not give a statement.
The methods developed in \cite{Sch02,Sch04} serve to check the conditions in Theorem~\ref{14} algorithmically.
We have used the {\sc Singular} \cite{GPS05} implementation \cite{Sch04a} of these methods to compute the following example which is out of the scope of Theorem~\ref{19}.

\begin{exa}
Consider the isolated singularity $D$ defined by $f=x^5+x^2y^2+y^5+z^5$.
By a Gr\"obner basis computation, one easily verifies that $f\not\in\ideal{\frac{\del f}{\del x},\frac{\del f}{\del z},\frac{\del f}{\del z}}$ which shows $D$ is not quasihomogeneous.
The spectrum of $f$ consists of the collection of $\alpha\in\QQ$ with multiplicity $\mu_\alpha\in\NN$ listed in Table~\ref{21}.
As there are no integer spectral numbers, the monodromy does not have an eigenvalue $1$. 
Thus, Theorem~\ref{14} implies that the logarithmic comparison theorem does not hold for $D$.
\begin{longtable}[t]{c|c|c|c|c|c|c|c|c|c|c|c|c|c}
\caption{Spectrum of $f=x^5+x^2y^2+y^5+z^5$}\label{21}\\
$\alpha$ & $-\frac{3}{10}$ & $-\frac{1}{10}$ & $\frac{1}{10}$ & $\frac{1}{5}$ & $\frac{3}{10}$ & $\frac{2}{5}$ & $\frac{1}{2}$ & $\frac{3}{5}$ & $\frac{7}{10}$ & $\frac{4}{5}$ & $\frac{9}{10}$ & $\frac{11}{10}$ & $\frac{13}{10}$ \\\hline
$\mu_\alpha$ & $1$ & $3$ & $5$ & $1$ & $7$ & $1$ & $8$ & $1$ & $7$ & $1$ & $5$ & $3$ & $1$
\end{longtable}
\end{exa}

We shall prove Theorem~\ref{14} in Section~\ref{18} after some preparations on logarithmic vector fields in the following Section~\ref{17}.

\section{Logarithmic vector fields}\label{17}

We shall assume throughout that $D$ is an isolated singularity and use the notation in \eqref{22}.
Denote by
\[
\Der(-\log D):=\{\delta\mid\delta(f)\in\calO f\}\subseteq\Der:=\Der_\CC(\calO)\cong\calO^{n+1}
\]
the $\calO$-module of logarithmic vector fields along $D$.
We may assume that
\[
\Der(-\log D)\subseteq\mm\Der=:\Delta
\]
which means that $D\not\cong D'\times\CC$.
Let $\delta_0$ be the image of $\delta\in\Delta$ and $\Der(-\log D)_0$ that of the infinite Lie algebra $\Der(-\log D)$ under the Lie algebra homomorphism
\begin{equation}\label{23}
\xymat{\Delta\ar[r]^-{\pi_0}&\Delta/\mm\Delta=:\Delta_0\cong\gl_{n+1}(\CC)},\quad\sum_{i,j=1}^na_{j,i}(x)x_i\del_j\mapsto(a_{i,j}(0))_{i,j},
\end{equation}
where we abbreviate $\del_i:=\frac{\del}{\del x_i}$ for $i=0,\dots,n$.
Note that $\Der(-\log D)_0$ is a finite Lie algebra.
The basis $x=x_0,\dots,x_n$ of $\mm$ defines a section of the map $\pi_0$ in \eqref{23} by which we can consider $\Delta_0$ and $\Der(-\log D)_0$ as Lie subalgebras of $\Delta$.
We call $\delta=\delta_0\in\Delta$ semisimple if the corresponding matrix $\pi_0(\delta)$ has this property.
If $\tau_0(\delta)$ is a nilpotent matrix (but not necessarily $\delta=\delta_0$), we call $\delta\in\Delta$ nilpotent.
While semisimplicity depends on the coordinate system, nilpotency is an intrinsic property.
Any $\delta\in\Delta$ can be decomposed as
\begin{equation}\label{24}
\delta=\delta_S+\delta_N
\end{equation}
into a semisimple $\delta_S$ and nilpotent $\delta_N$.
Note that $\delta_S$ is just the semisimple part of $\tau_0(\delta)$ in the linear algebra sense.

\begin{prp}\label{1}
If $D$ is of order at least $3$ (which means that $f\in\mm^3$) but not quasihomogeneous then $\Der(-\log D)$ contains only nilpotent vector fields.
\end{prp}

\begin{proof}
We may replace $\Der(-\log D)$ by its $\mm$-adic completion $\widehat\Der(-\log D)=\Der(-\log\widehat D)$ where $\widehat D$ is defined by the same $f\in\calO\subseteq\widehat\calO$ considered as a formal power series.
Indeed, $D$ has an isolated singularity if and only if $\widehat D$ has and quasihomogeneity is equivalent to Euler homogeneity by \cite{Sai71} which is invariant under completion.
Moreover, $\widehat\Der(-\log D)_0=\Der(-\log D)_0$ and the notion of nilpotency is preserved.

Let $\widehat\delta\in\Der(-\log\widehat D)$ and decompose it as in \eqref{24}.
By \cite[Thm.~5.4]{GS06}, there is a (formal) coordinate system with respect to which $\sigma:=\widehat\delta_S\in\Der(-\log\widehat D)$ and a defining equation $\widehat f\in\widehat\calO$ of $\widehat D$ such that $\sigma(\widehat f)\in\QQ\widehat f$.
We must have $\sigma(\widehat f)=0$ as otherwise $\widehat D$ would be quasihomogeneous by \cite{Sai71}. 
Assume that $\sigma\ne0$.
This means that the monomial support of $\widehat f$ lies in a proper vector subspace.

As $\widehat f$ has an isolated critical point, \cite[Cor.~1.6]{Sai71} states that, for each $j=0,\dots,n$, there must be a monomial with exponent $me_j$ or $me_j+e_{j'}$ in the monomial support of $\widehat f$.
But, by the order hypothesis, $\widehat f\in\widehat \mm^3$ which implies that these monomials are linearly independent.
This contradicts to the monomial support of $\widehat f$ having codimension at least one and finishes the proof.
\end{proof}

Dropping the order hypothesis in Proposition~\ref{1}, a weaker statement holds.

\begin{prp}\label{15}
If $D$ is not quasihomogeneous then $\tr(\delta_0)=0$ for all $\delta\in\Der(-\log D)$.
\end{prp}

\begin{proof}
If $f\in\mm^3$ then we may assume by Proposition~\ref{1} that $\delta_0$ is a lower triangular matrix and the claim follows.

In the general case, we can assume by the Splitting Lemma that
\[
f(x)=f'(x')+\sum_i(x''_i)^2,\quad x=(x',x''),\quad\del=(\del',\del''),\quad f'\in\mm^3.
\]
With $D$ also
\[
D'=\{x'\mid f'(x')=0\}\subseteq X'
\]
is a nonquasihomogeneous isolated singularity by \cite{Sai71}.
Writing $\delta=\sum_ig_i\del_i$, we have to check that the monomial $x_i$ does not occur in $g_i$.
By definition of $\Der(-\log D)$, $\delta$ corresponds to a syzygy of
\begin{gather}
\label{16a}x'_i\del'_j(f)=x'_i\del'_j(f')\in\mm^3,\\
\label{16b}x''_i\del''_j(f)=2x''_ix''_j\in\mm^2\smallsetminus\mm^3,\\
\label{16c}x''_i\del'_j(f)=x''_i\del'_j(f')\in\mm^3,\quad
x'_i\del''_j(f)=2x'_ix''_j\in\mm^2\smallsetminus\mm^3,\quad
x'_if,x''_jf\in\mm^3.
\end{gather}
By \cite{Sai71}, $f$ can not occur with a constant coefficient in \eqref{16c} as $D$ is assumed not to be quasihomogeneous.
We are concerned only with the constant coefficients of \eqref{16a} and \eqref{16b} for $i=j$.
Those of \eqref{16b} are obviously zero.
Setting $x''=0$ yields a syzygy of \eqref{16a} and $x'_if'$ that induces an element of $\Der(-\log D')$.
Thus, the constant coefficients of \eqref{16a} are zero for $i=j$ by the first part of the proof.
\end{proof}

Let $\Omega_X^\bullet$ be the complex of holomorphic differential forms on $X$ and denote the volume form by
\[
dx:=dx_0\wedge\dots\wedge dx_n\in\Omega_X^{n+1}.
\]
The complex of logarithmic differential forms along $D$ was introduced in \cite{Sai80} as
\begin{equation}
\Omega^\bullet(\log D):=\Omega^\bullet_X(D)\cap d^{-1}\Omega^\bullet_X(D)\subseteq\Omega_X^\bullet(*D).
\end{equation}

\begin{cor}\label{2}
If $D$ is not quasihomogeneous then $0\ne[\frac{udx}{f}]\in\calH^{n+1}(\Omega^\bullet(\log D))$ for any $u\in\calO^*$.
\end{cor}

\begin{proof}
The module $\Omega^n(\log D)$ is the image of the inner product 
\[
\xymat@R=0pt{
\Der(-\log D)\times\Omega^{n+1}(\log D)\ar[r]&\Omega^n(\log D)\\
(\delta,\omega)\ar@{|->}[r]&\iota_\delta(\omega)
}
\]
defined in \cite[Lem.~1.6.ii]{Sai80}.
As $\Omega^{n+1}(\log D)=\calO_X(D)dx$, we have $\Omega^n(\log D)=\iota_{\Der(-\log D)}dx/f$.
Let $\delta=\sum_ig_i\del_i\in\Der(-\log D)$ and note that $\delta(f)\in\mm f$ by nonquasihomogeneity of $D$ and \cite{Sai71}.
Then we compute
\begin{align*}
fd(\iota_\delta(dx)/f)&=d\circ\iota_\delta(dx)-(df/f)\wedge\iota_\delta(dx)\\
&=\Lie_\delta(dx)-(\delta(f)/f)dx\\
&=\sum_i\del_i(g_i)dx-(\delta(f)/f)dx\equiv\tr(\delta_0)dx\mod\mm.
\end{align*}
By Proposition~\ref{15}, this implies that $d\Omega^n(\log D)\subseteq\mm\Omega^{n+1}(\log D)=\mm\Omega^{n+1}_X(D)$ and the claim follows.
\end{proof}

\section{Gauss--Manin system}\label{18}

We keep our general assumption that $D$ is an isolated singularity and continue to use the notation in \eqref{22}.
Corollary~\ref{2} leads us to study the necessary condition
\begin{equation}\label{9}
0\ne\left[\frac{udx}{f}\right]\in\calH^{n+1}(\Omega_X^\bullet(*D))\text{ for all }u\in\calO^*
\end{equation}
for the logarithmic comparison theorem to hold for nonquasihomogeneous $D$.
We shall reformulate this condition in terms of the Gauss--Manin system of $f\colon X\to T$ using \cite[\S1-2]{Kar91} as a starting point.

Let $M$ be the monodromy on the canonical Milnor fiber $X_\infty$ of $f$ \cite[\S5]{SS85}.
By construction of $X_\infty$ and \cite{Mil68},
\begin{equation}\label{27}
H^k(X_\infty;\CC)\cong H^k(X_t;\CC)=0,\quad\text{ if }k\ne0,n,
\end{equation}
where $X_t:=f^{-1}(t)$ and $t\in T^*:=T\backslash\{0\}$.
Then the cohomological Wang sequence reads
\begin{equation}\label{4}
\xymat@C=16pt{
0\ar[r]&H^n(U;\CC)\ar[r]&H^n(X_\infty;\CC)\ar[r]^-{M-1}&H^n(X_\infty;\CC)\ar[r]&H^{n+1}(U;\CC)\ar[r]&0.}
\end{equation}
Recall that the eigenvalues of $M$ on $H^n(X_\infty;\CC)$ are roots of unity by the monodromy theorem \cite[Satz 4]{Bri70}. 
Decompose $M=M_sM_u$ into semisimple and unipotent part and let $H^k(X_\infty;\CC)_\rho$ denote the generalized $\rho$-eigenspace of $M$. 
Then $M-1$ has the same kernel and cokernel on $H^n(X_\infty;\CC)$ as on $H^n(X_\infty;\CC)_1$, $M$ coincides with $M_u$ on $H^n(X_\infty;\CC)_1$, and $M_u-1$ has the same kernel and cokernel as $N:=\log M_u$ on $H^n(X_\infty;\CC)_1$.
Thus, \eqref{4} leads to an exact sequence
\begin{equation}\label{5}
\xymat@C=10pt{
0\ar[r]&H^n(U;\CC)\ar[r]&H^{n+1}(X_\infty;\CC)_1\ar[r]^-N&H^{n+1}(X_\infty;\CC)_1\ar[r]&H^{n+1}(U;\CC)\ar[r]&0}.
\end{equation}

To see the D-module structure hidden in \eqref{5} requires a refined approach.
Let $\Gamma$ be the graph of $f$ and consider the maps $i(x)=(x,0)$, $j(x)=(x,f(x))$, and $p(x,t)=t$ in the diagram
\[
\xymat{
X\ar[r]^-j&\Gamma\ar@{^(->}[r]&X\times T\ar[d]^-p&X\ar[l]_-i\\
&&T.
}
\]
Then $i(X)\cap\Gamma=D$, $i^*\Omega_{X\times T/T}^\bullet(*\Gamma)=\Omega_X^\bullet(*D)$, and there is an exact sequence
\begin{equation}\label{28}
\xymat{
0\ar[r]&\Omega_{X\times T/T}^\bullet(*\Gamma)\ar[r]^-t&\Omega_{X\times T/T}^\bullet(*\Gamma)\ar[r]&i_*\Omega_X^\bullet(*D)\ar[r]&0}.
\end{equation}
As $X$ is Stein and $\Omega^\bullet_{X\times T/T}$ consists of $\calO_{X\times T}$-coherent and hence $p_*$-acyclic modules, the Poincar\'e Lemma shows that
\[
\RR^kp_*\Omega^\bullet_{X\times T/T}=0,\quad\text{ if }k\ge1.
\]
Therefore
\begin{align*}
\RR^kp_*(\Omega^\bullet_{X\times T/T}(*\Gamma))
&=\RR^kp_*(\Omega^\bullet_{X\times T/T}(*\Gamma)/\Omega^\bullet_{X\times T/T})\\
&=\RR^kf_*j^{-1}(\Omega^\bullet_{X\times T/T}(*\Gamma)/\Omega^\bullet_{X\times T/T})=\int_f^{k-(n+1)}\calO_X,\quad\text{ if }k\ge1,
\end{align*}
is the Gauss-Manin system of $f\colon X\to T$.
As $\int_f^{k-n}\calO_X$ is a $\calD_T$-coherent regular extension of $\calO_{T^*}\left(\bigcup_{t\in T^*}H^k(X_t;\CC)\right)$ to $T$, it follows from \eqref{27} that
\[
\int_f^k\calO_X=0,\quad\text{ if }k\ne 0,-n.
\]
Using that $\RR p_*i_*=\RR p_*\RR i_*=\RR(p\circ i)_*=\RR0_*=\RR\Gamma(X,\_)=\HH(X,\_)$, Grothendieck's comparison theorem \cite{Gro66} implies that
\[
\RR^kp_*i_*\Omega_X^\bullet(*D)=\HH^k(X,\Omega_X^\bullet(*D))=h^k(\Gamma(X,\Omega_X^\bullet(*D)))=\calH^k(\Omega_X^\bullet(*D)).
\]
So applying $\RR p_*$ to \eqref{28} yields an exact sequence  
\begin{equation}\label{3}
\xymat@R=0pt{
0\ar[r]&\calH^n(\Omega_X^\bullet(*D))\ar[r]&\calG\ar[r]^-t&\calG\ar[r]^-\pi&\calH^{n+1}(\Omega_X^\bullet(*D))\ar[r]&0\\
&&&[udx]\ar@{|->}[r]&[\frac{udx}{f}]}
\end{equation}
where $\calG:=\int_f^0\calO_X$.

By \cite[\S15]{Pha79} and \cite[Lem.~3.3]{SS85}, $\calG$ can be represented in explicit form as
\begin{gather}\label{11}
\calG\cong\Omega_X^{n+1}[D]/(d-D\cdot df\wedge)\Omega_X^n[D],\\
\nonumber
\del_t\left[\frac{\omega}{(f-t)^k}\right]=\left[\frac{k!\omega}{(f-t)^{k+1}}\right]\mapsto[\omega D^k].
\end{gather}
The operator $\del_t$ is invertible on $\calG$ by \cite[\S15.2.2]{Pha79} and $\calG\cong\Cs[\del_t]^\mu$ where $\mu$ is the Milnor number of $f$ and $\Cs$ is the ring of microdifferential operators with constant coefficients.
Composing the second map in \eqref{3} with $\del_t^{-1}$, the operator $t$ in \eqref{3} can be replaced by $t\del_t$ without changing the cokernel.
This leads to an exact sequence 
\begin{equation}\label{6}
\xymat{
0\ar[r]&\calH^n(\Omega_X^\bullet(*D))\ar[r]&C^0\ar[r]^-{t\del_t}&C^0\ar[r]^-\pi&\calH^{n+1}(\Omega_X^\bullet(*D))\ar[r]&0}
\end{equation}
where $C^\alpha$ denotes the generalized $\alpha$-eigenspace of the operator $t\del_t$ on $\calG$.
The (decreasing) Kashiwara--Malgrange filtration $V^\bullet$ on $\calG$ is essentially defined by $\gr_V^\alpha\calG\cong C^\alpha$ and consists of free $\Cs$-modules of rank $\mu$.
By \cite[\S5, p.\,652]{SS85}, one can identify the vector spaces with endomorphisms
\begin{equation}\label{26}
(H^n(X_\infty;\CC)_\lambda,N)\cong(C^\alpha,-2\pi i(t\del_t-\alpha)).
\end{equation}
We can thus identify $N=-2\pi it\del_t$ in the sequences \eqref{5} and \eqref{6}.
By \cite{Gro66}, also the outer terms of these sequences coincide.

With \eqref{9} and \eqref{3} in mind, we are interested in the image of the canonical map $\Omega_X^{n+1}\to\calG$ (see \eqref{11}), which is the Brieskorn lattice 
\[
\calH''\cong\Omega_X^{n+1}/df\wedge d\Omega_X^{n-1}
\]
of $f$ \cite{Bri70}.
From \eqref{15} it follows easily that
\begin{equation}\label{30}
\calH''/\del_t^{-1}\calH''\cong\Omega_X^{n+1}/df\wedge\Omega_X^n=:\Omega_f\cong\CC^\mu.
\end{equation}
By \cite{Seb70}, $\calH''$ is a free $\CC\{t\}$-module of rank $\mu$ and, by \cite[Lem.~4.5]{Mal74},
\begin{equation}\label{29}
\calH''\subseteq V^{>-1}
\end{equation}
from which one can derive that $\calH''$ is also a free $\Cs$-module of rank $\mu$ \cite[Prop.~2.5]{Pha77}.

For $g\in\calG$ and $\alpha\in\QQ$, we shall write $g_\alpha$ for the $C^\alpha$-component of $g$.
The preceding arguments now show that \eqref{9} is equivalent to
\begin{equation}\label{10}
[udx]_0\not\in N(C^0)\text{ for all }u\in\calO^*.
\end{equation}
By \eqref{17}, condition \eqref{14a} in Theorem~\ref{14} implies that $C^0=0$ and the claim follows in that case.
The spectrum of $f$ is defined as the spectrum $\alpha_1\le\cdots\le\alpha_\mu$ of the filtration induced by $V^\bullet$ on $\Omega_f$, that is,
\[
\#\{i\mid\alpha=\alpha_i\}=\dim_\CC\gr_V^\alpha\Omega_f.
\]
Under condition \eqref{14b} in Theorem~\ref{14}, $\calH''\subseteq V^{>0}$ and hence $[udx]_0=0$ for all $u\in\calO$.
Thus, also in this case, Theorem~\eqref{14} holds true.

In order to prove Theorem~\ref{14} under the assumption \eqref{14c} or \eqref{14d}, let us assume that $\alpha_1<0$.
From \eqref{30} and \eqref{29}, we conclude that $C^0\cap\del_t^{-1}\calH''=0$ and hence
\begin{equation}\label{13}
C^0\subseteq V^{>-1}/\del_t^{-1}\calH''\supseteq\calH''/\del_t^{-1}\calH''\cong\Omega_f.
\end{equation}
By \cite[Rem.~3.11]{Sai91}, $\mm dx$ surjects onto $V^{>\alpha_1}\Omega_f$.
In particular, $\alpha_1<\alpha_2$ and $[udx]\in V^{\alpha_1}\setminus V^{>\alpha_1}$ for all $u\in\calO^*$.
Moreover, $\gr_V^0\calH''=:H_0\subseteq C^0$ is in the image of $\mm dx$ by \eqref{13}.
By \cite[Lem.~3.4 and \S6.5]{SS85}, $\del_t^{n-k}\calH''$ induces the Hodge filtration $F^k$ on $\gr_V^0\calG=C^0$ for which $N$ is a morphism of type $-1$.
Therefore, $H_0$ has a complement $G_0$ in $C^0$ such that $N(H_0)\subseteq G_0$.
This shows that \eqref{10} is equivalent to 
\begin{equation}\label{12}
[udx]_0\not\in H_0\oplus N(C^0)\text{ for all }u\in\calO^*.
\end{equation}
This proves Theorem~\ref{14} under the hypothesis \eqref{14c}.

Assume finally that $\alpha_2=0$.
Then, modulo $\CC^*$, the $G_0$-component of $[udx]_0$ is independent of $u\in\calO^*$ and \eqref{12} is equivalent to
\begin{equation}\label{16}
[dx]_0\not\in H_0\oplus N(C^0).
\end{equation}
This finishes the proof of our main result Theorem~\ref{14}.

\bibliographystyle{amsalpha}
\bibliography{lctis}

\end{document}